# On Equal Point Separation by Planar Cell Decompositions

Nikhil Marda

**Abstract**

In this paper, we investigate the problem of separating a set $X$ of points in $\mathbb{R}^2$ with an arrangement of $K$ lines such that each cell contains an asymptotically equal number of points (up to a constant ratio). We consider a property of curves called the stabbing number, defined to be the maximum countable number of intersections possible between the curve and a line in the plane. We show that large subsets of $X$ lying on Jordan curves of low stabbing number are an obstacle to equal separation. We further discuss Jordan curves of minimal stabbing number containing $X$. Our results generalize recent bounds on the Erdős-Szekeres Conjecture, showing that for fixed $d$ and sufficiently large $n$, if $|X| \geq 2^{c_d n/d + o(n)}$ with $c_d = 1 + O(\frac{1}{\sqrt{d}})$, then there exists a subset of $n$ points lying on a Jordan curve with stabbing number at most $d$.

# 1 Introduction

A classic result in combinatorial geometry is the Szemerédi-Trotter theorem, which bounds the number of incidences between a set $X$ of points and a set $L$ of lines in the plane. An elegant proof [4] is based on the method of cell decomposition. This involves placing an arrangement $\mathcal{A}_K$ of $K$ lines in the plane such that each cell created by $\mathcal{A}_K$ has roughly the same number of incidences with either $X$ or $L$. Then the number of incidences between $X$ and $L$ within each cell can be estimated and multiplied by the number of cells, giving the desired bound.

Guth [8] notes that while $\mathcal{A}_K$ can always be found such that each cell has roughly the same number of incidences with $L$, equal separation cannot always be done with $X$. He asks when it is possible to find $\mathcal{A}_K$ such that each cell contains roughly the same number of points in $X$, such as in Figure 1.1.

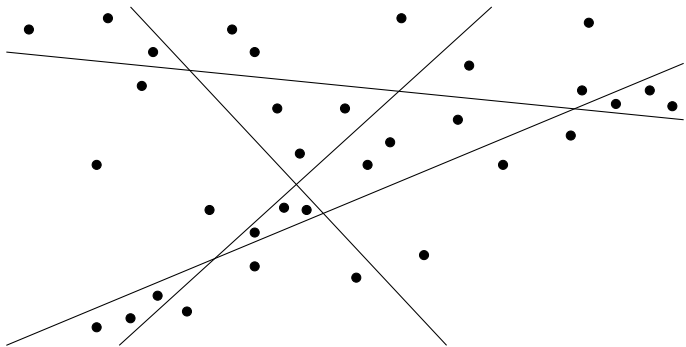

**Figure 1.1:** Example of $X$ separated equally by the cells of $\mathcal{A}_4$.

Hereafter, we let $X$ denote a finite set of $N$ points in general position (no 3 points are collinear) and we let $\mathcal{A}_K$ denote a simple arrangement of $K$ lines (no 2 lines are parallel and no 3 lines meet at a point, thus creating $\frac{K^2+K+2}{2}$ cells). We say $A(N,K) \sim B(N,K)$ if there exist constants $M$, $C_1$, $C_2 > 0$ such that when $N, K \geq M$, then $C_1 B(N,K) \leq A(N,K) \leq C_2 B(N,K)$. We also say $A(N,K) \ll B(N,K)$ if for all $C > 0$, there exists $M$ such that when $N, K \geq M$, we have $|A(N,K)| < C|B(N,K)|$. We will also use $A(N,K) \lesssim B(N,K)$ to imply that either $A(N,K) \ll B(N,K)$ or $A(N,K) \sim B(N,K)$. We are exploring equal separation by $\mathcal{A}_K$: Given some $X$, the lines of $\mathcal{A}_K$ split the plane such that each of the $\sim K^2$ cells contains $\sim \frac{N}{K^2}$ points.

Guth [8] mentions this problem is nontrivial because there are $\sim K$ parameters but $\sim K^2$ conditions on the cells to satisfy. He mentions an explicit counterexample to equal separation. Consider when $X$ lies on a convex curve $J$, such as in Figure 1.2. Every line $\ell \in \mathcal{A}_K$ intersects $J$ at most twice. Hence, the $K$ lines split $J$ into at most $2K$ segments, each of which is contained entirely in some cell. Thus, there exists a cell with at least $\frac{N}{2K}$ points, which is greater than $\frac{N}{K^2}$.

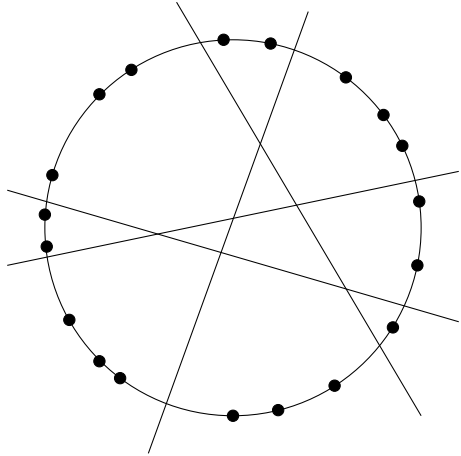

**Figure 1.2:** Example of $\mathcal{A}_4$ and $X$ on a convex curve $J$.

In this paper, we search for obstacles to equal separation. In the process, we build intuition about configurations preventing equal separation. We approach the problem by first introducing the following quantity.

**Definition 1.1 (Cutting number).** The cutting number of $X$ is the minimum number of points that must lie in the cell containing the most points for any $\mathcal{A}_K$:

$$\mathrm{CUT}_\mathrm{K}(X) = \min_{\mathcal{A}_K}\left(\max_{\text{cells}} |X \bigcap \text{cell}|\right).$$

It was noted earlier that if $X$ lies on a convex curve, then $\mathrm{CUT}_\mathrm{K}(X) \geq \frac{N}{2K} \gg \frac{N}{K^2}$. Also by the definition of equal separation, for any $X$, if $\mathrm{CUT}_\mathrm{K}(X) \gg \frac{N}{K^2}$, then $X$ cannot be separated equally.

For $X$ on a convex curve $J$, we used the fact that any line intersects $J$ at most 2 times to show that $\mathrm{CUT}_\mathrm{K}(X)$ was large. This can be generalized by considering curves that have a small number of intersections with lines in the plane.

**Definition 1.2 (Stabbing number).** For an object $J$, the stabbing number $\mathrm{stab}(J)$ is the maximum countable number of intersections possible between the curve and a line in the plane.

The terms "crossing numbers" and "stabbing numbers" both appear in the literature, but to avoid confusion with crossing numbers from graph theory, we use the latter terminology. Stabbing numbers have been studied in other contexts and have applications in discrepancy theory [13], range searching [3], approximating zonoids with zonotopes [10], and finding minimum cost spanning trees [1]. As a consequence of their relation with our original problem, we investigate low stabbing number Jordan curves (plane simple closed curves) containing $X$. These results build toward solving a problem in Ramsey theory: a generalization of recent bounds [15] on the Erdős-Szekeres Conjecture. We show that for fixed $d$ and sufficiently large $n$, given $2^{c_d n/d + o(n)}$ points with $c_d = 1 + O(\frac{1}{\sqrt{d}})$, there exists a subset of $n$ points lying on a Jordan curve with stabbing number at most $d$.

Section 2 of this paper addresses upper bounds on the cutting number. Section 3 demonstrates the relevance of low stabbing number Jordan curves to lower bounds on the cutting number. Section 4 summarizes the obstacles discovered to equal separation and expands upon our discussion of minimal stabbing number Jordan curves. Section 5 introduces the Erdős-Szekeres Conjecture and generalizes it to Jordan curves of bounded stabbing number. Finally, Section 6 concludes with thoughts for future research and the relevance of our work in the literature.

## 2 Upper bounds on the cutting number

We first seek to understand the behavior of the cutting number. As noted earlier, if $X$ lies on a convex curve, then its cutting number is at least $\frac{N}{2K}$. In this section, we show that $X$ lying on a convex curve actually gives the largest possible cutting number. This is done by proving the existence, for all $X$, of a line arrangement $\mathcal{A}_K$ such that no cell has more than $\lceil \frac{N}{2K} \rceil$ points. The following definition and lemma will be useful.

**Definition 2.1 (Linearly separable).** We call two sets of points $R$ and $S$ linearly separable if there exists a line such that all points in $R$ lie on one side of the line and all points in $S$ lie on the other side of the line.

**Lemma 2.1.** *Let $A$ and $B$ be linearly separable sets of at least $r$ points each. Then there exists a line which forms a half-plane containing exactly $r$ points of $A$ and $r$ points of $B$.*

*Proof.* Let $X = A \cup B$ and let $A$ and $B$ be separated by a line $j$. Consider a unit vector $\vec{u}$ rotating counterclockwise, skipping over the finitely many $\vec{u}$ that are perpendicular to a line containing 2 points in $X$. Let $\ell$ be a line perpendicular to $\vec{u}$ such that the open half-plane $\mathcal{P}_\ell$ in the direction of $\vec{u}$ contains $2r$ points. Define $f(\ell)$ to be the number of points of $A$ in $\mathcal{P}_\ell$ for some $\ell$. When $\vec{u}$ is perpendicular to $j$ and directed toward the half-plane containing $A$, we have $f(\ell) \geq r$.

As we rotate $\vec{u}$ counterclockwise, $f(\ell)$ can only change when the rotating $\vec{u}$ passes over a direction perpendicular to a line through 2 elements of $X$. When this happens, some $x_i \in X$ may be removed from $\mathcal{P}_\ell$ and some $x_j \in X$ may be added to $\mathcal{P}_\ell$. This can change $f(\ell)$ by at most 1. When $\vec{u}$ is perpendicular to $j$ and directed toward the half-plane containing $B$, we have $f(\ell) \leq r$. It follows that there exists an $\ell$ for which $f(\ell) = r$. □

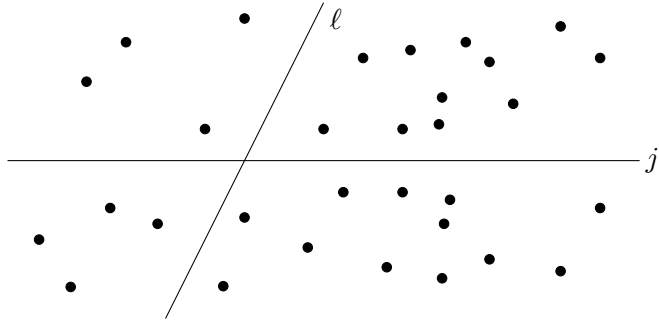

**Figure 2.1:** Example of Lemma 2.1 applied on $A$ and $B$ for $r = 4$.

An example is shown in Figure 2.1. Given enough lines, we can ensure no cell has more than a certain number of points.

**Theorem 2.2.** *If $K \geq \left\lceil \frac{N}{2H} \right\rceil$, then $\mathrm{CUT}_K(X) \leq H$.*

*Proof.* We will construct an $\mathcal{A}_K = \{\ell_1, \ell_2, ..., \ell_K\}$ iteratively such that no cell has more than $H$ points. Let $\ell_1$ be a line disjoint from $X$ separating $X$ into 2 cells, each with at most $\left\lceil \frac{N}{2} \right\rceil$ points. For each side of $\ell_1$, consider the set of points which lie in a cell with more than $H$ points for $\mathcal{A}_i = \{\ell_1, \ell_2, ..., \ell_i\}$, and denote these sets $R_i$ and $S_i$.

Using Lemma 2.1, we take $\ell_{i+1}$ for $i \geq 1$ to separate off 1 cell of $H$ points from both $R_i$ and $S_i$. It follows $|R_i|$ and $|S_i|$ are at most $\left\lceil \frac{N}{2} \right\rceil - (i-1)H$. We want $\left\lceil \frac{N}{2} \right\rceil - (i-1)H \leq H$. Rearranging, we get $i \geq \frac{\lceil N/2 \rceil}{H}$. Since $i$ is an integer, we have $i \geq \left\lceil \frac{\lceil N/2 \rceil}{H} \right\rceil$ and it can be shown $\left\lceil \frac{\lceil N/2 \rceil}{H} \right\rceil = \left\lceil \frac{N}{2H} \right\rceil$. □

This proof also gives us a way to form subsets of $X$ with equal magnitude and disjoint convex hulls. When every line $\ell_i$ is added, it forms 2 cells with $H$ points. However, $\ell_i$ might intersect cells of $H$ points formed by $\ell_j$ for $j < i$. Since these cells already have $H$ points, we can delete the parts of the new lines in these cells. Then we will have at least $2H$ cells with exactly $H$ points each.

**Definition 2.2 (Partial cutting).** For $\ell_i \in \mathcal{A}_K$ with $i \geq 2$ constructed in the proof of Theorem 2.2, let the half-plane containing $2H$ points $\ell_i$ separates be $Q_i$. For all $s < r$, delete portions of $\ell_r$ intersecting $Q_s$. We call this construction a partial cutting.

An example of a partial cutting is shown in Figure 2.2. Every $\ell_i$ for $i \geq 3$ becomes either a ray or line segment. Partial cuttings will be helpful in Section 5.

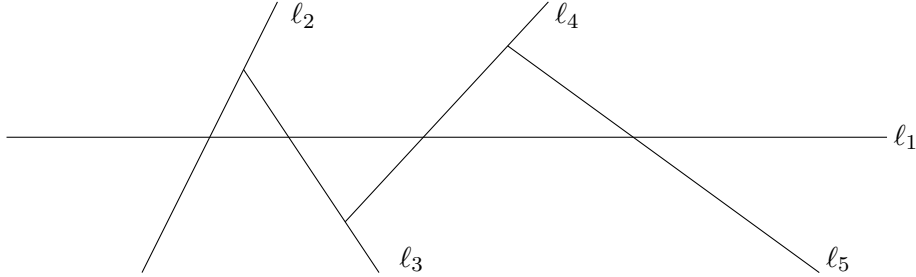

**Figure 2.2:** Example of a partial cutting.

Our bound on the cutting number from Theorem 2.2 can be rephrased for $X$ in terms of $N$ and $K$. This is useful because we want to compare $\text{CUT}_K(X)$ to $\frac{N}{K^2}$.

**Corollary 2.2.1.** *For any $X$, we have $\text{CUT}_K(X) \leq \left\lceil \frac{N}{2K} \right\rceil$.*

*Proof.* Take $H = \left\lceil \frac{N}{2K} \right\rceil$. It can be shown that $\left\lceil \frac{N}{2\lceil N/(2K)\rceil} \right\rceil = K \geq K$, so by Theorem 2.2, it follows that $\text{CUT}_K(X) \leq \left\lceil \frac{N}{2K} \right\rceil$. $\square$

This is a tight upper bound on the cutting number.

**Corollary 2.2.2.** *If $X$ lies on a convex curve, then $\text{CUT}_K(X) = \left\lceil \frac{N}{2K} \right\rceil$ for all $K$.*

*Proof.* We showed in Section 1 that if $X$ lies on a convex curve, then $\text{CUT}_K(X) \geq \frac{N}{2K}$. Since $\text{CUT}_K(X)$ is an integer, it follows from Corollary 2.2.1 that $\text{CUT}_K(X) = \left\lceil \frac{N}{2K} \right\rceil$. $\square$

Finally, we want to make a statement asymptotically. That is, if a "large" subset of $X$ lies on a convex curve, then can we say something about the cutting number? We first define the following.

**Definition 2.3 (Convex number).** The convex number of $X$ is the maximum number of points that lie on some convex curve in the plane:
$$\text{CON}(X) = \max_{\text{convex curves}} |X \cap \text{curve}|.$$

**Corollary 2.2.3.** *If $\text{CON}(X) \sim N$, then $\text{CUT}_K(X) \sim \frac{N}{K}$.*

*Proof.* By Corollary 2.2.1, we have $\text{CUT}_K(X) \leq \left\lceil \frac{N}{2K} \right\rceil \lesssim \frac{N}{K}$. If $\text{CON}(X) \sim N$, then there exists $Y \subseteq X$ such that $\text{CUT}_K(X) \geq \text{CUT}_K(Y) \gtrsim \frac{N}{K}$. This gives $\text{CUT}_K(X) \sim \frac{N}{K}$. $\square$

A large subset of $X$ lying on a convex curve is an obstacle to equal separation. In the next section, we show that it is not the only such obstacle.

## 3 Low stabbing number Jordan curves

In Section 1, we suggested that $X$ lying on a Jordan curve with low stabbing number was also an obstacle to equal separation. In this section, we flesh out this idea. Furthermore, we show a low stabbing number Jordan curve containing $X$ implies a large value of $\text{CUT}_\text{K}(X)$ but does not necessarily mean $\text{CON}(X)$ is large. That is, there exists an obstacle to equal separation distinct from the presence of a large convex subset.

**Definition 3.1 ($d$-curve).** A $d$-curve is a Jordan curve $J$ with $\text{stab}(J) \leq d$.

For example, convex curves (2-curves) are also 4-curves. The relationship between the cutting number and $d$-curves containing $X$ follows as expected.

**Lemma 3.1.** *If a $d$-curve $Y$ contains $X$, then $\text{CUT}_\text{K}(X) \geq \frac{N}{Kd}$.*

*Proof.* Any line intersects $Y$ at most $d$ times, so $|\mathcal{A}_K \cap Y| \leq Kd$. At most $Kd$ segments are formed along $Y$, each contained in a cell. Thus some cell has $\geq \frac{N}{Kd}$ points. □

Is the presence of a large convex subset the only obstacle to equal separation? Corollary 2.2.3 told us if $\text{CON}(X) \sim N$, then $\text{CUT}_\text{K}(X) \sim \frac{N}{K} \gg \frac{N}{K^2}$. However, we can show that there exists $X$ with large $\text{CUT}_\text{K}(X)$ and small $\text{CON}(X)$ by constructing an $X$ that lies on a curve with low stabbing number and contains no large convex subset.

**Theorem 3.2.** *There exists $X$ with $\text{CON}(X) \lesssim \sqrt[3]{N}$ and $\text{CUT}_\text{K}(X) \gtrsim \frac{N}{K}$.*

*Proof.* We first define some notation. In the $xy$-plane, let $C(r, a, b)$ denote the top half of a circle with radius $r$ centered at $(a, b)$. We define $C_0 := C(s, 0, 0)$ for some $s$ to be chosen later. Let $M$ be an integer such that $M \sim \sqrt[3]{N}$. For all nonzero integers $k$ such that $-M \leq k \leq M$, define $A_k := (\frac{s}{p})(\frac{k}{M+1})$ for some $p > 1$ to be chosen later. Then, for some $t \ll s$ to be chosen later, define

$$C_k := C(t, A_k, \sqrt{s^2 - A_k^2} - 2t).$$

Finally, having drawn $C_{-M}, ..., C_M$, we add vertical lines between $C_0$ and the endpoints of $C_k$, deleting the subtended parts of $C_0$. Let this final curve be called $G$. Figure 3.1 shows an example.

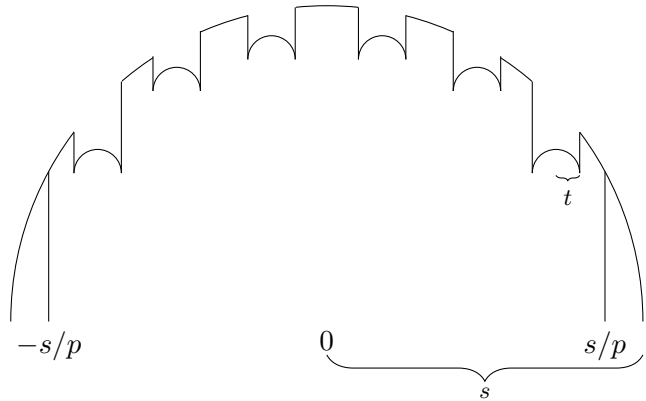

**Figure 3.1:** Example of $G$.

*Step 1.* We begin by bounding the maximum stabbing number of $G$ for sufficiently small $t$. Let $H_k$ be a circle of minimal radius such that it encloses $C_k$ and the vertical lines joining $C_k$ to $C_0$. We pick sufficiently small $t$ such that no line intersects more than two $H_k$. Then any line intersects $G$ at most 4 times in one $H_k$ and intersects $G$ at most 2 times outside of an $H_k$. This gives a total of at most 10 intersections between any line in the plane and $G$, making $G$ a 10-curve.

*Step 2.* Next, we place points on $G$. Suppose a convex curve contains three points on $C_k$, labeled $R, S, T$ by increasing $x$-coordinate. Then a fourth point $U$ cannot lie in the cell created by $\overrightarrow{RS}$ and $\overrightarrow{TS}$. An example is shown in Figure 3.2. On each $C_k$, we place $M$ points arbitrarily close to the intersection of $C_k$ with its tangent line parallel to $y = 0$. This means no point can be added from the half-plane opposite the tangent line parallel to $y = 0$. Thus for the points on $G$, a convex curve has $\lesssim M$ points from 2 $C_k$ and at most 2 points from $\lesssim M$ $C_k$, giving a convex curve with $\lesssim M$ points.

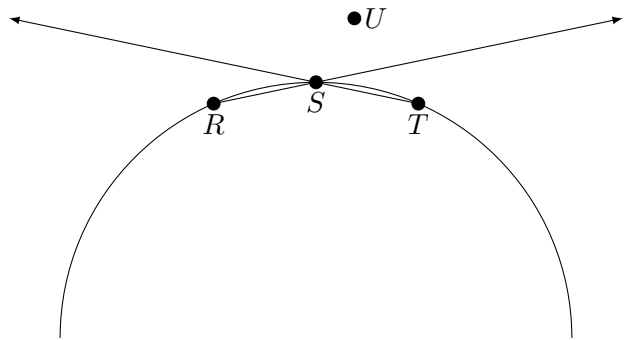

**Figure 3.2:** $R$, $S$, $T$, and $U$ cannot lie on the same convex curve. As $R$ and $T$ approach $S$, the region where $U$ cannot lie becomes the half-plane opposite the tangent line at $S$.

*Step 3.* Now we create a new curve of bounded maximum stabbing number. Consider the part of $G$ that lies above the line $y = \sqrt{s^2 - A_M^2} - 3t$. We place a copy $G_k$ of this part of the curve on each edge $E_k$ of a regular $M$-gon directed inward, where we choose $s$ to be small enough that no line in the plane intersects more than two $G_k$. We also delete the portion of $E_k$ subtended by $G_k$. Thus, any line in the plane has at most 10 intersections with each of two $G_k$ and at most 2 intersections with the $M$-gon, making this new curve a 22-curve. Figure 3.3 shows an example of this curve.

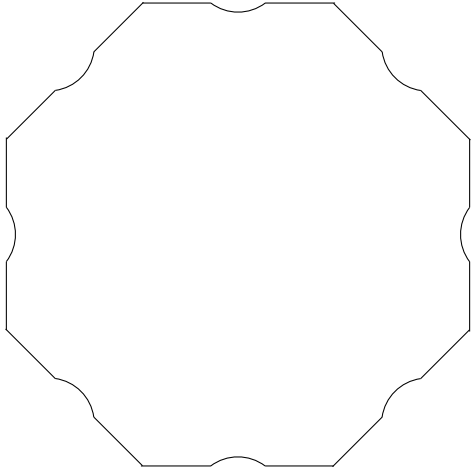

**Figure 3.3:** Example of discussed 22-curve.

*Step 4.* Finally, we show that we can pick $p$ such that this construction satisfies the statement. As $p$ goes to infinity, the tangent lines to $G_k$ approach the line parallel to $E_k$. Hence, for large enough $p$, no tangent line to $G_k$ intersects another $G_l$. Thus by the argument used in Step 2, if a convex curve contains at least 3 points from some $G_k$, it cannot contain points from any other $G_l$. This curve would contain $\lesssim M$ points as described earlier. Otherwise, at most 2 points can be taken from each $G_k$, which would also give a convex curve with $\lesssim M$ points.

Thus, we have $X$ with $M^3 \sim N$, points, with $\text{CON}(X) \lesssim \sqrt[3]{N}$ and $\text{CUT}_K(X) \geq \frac{N}{22K} \gtrsim \frac{N}{K}$.
$\square$

This shows that a large convex subset of $X$ is not the only obstacle to equal separation. Instead, we should look at $d$-curves containing $X$.

# 4 Minimal stabbing number Jordan curve through $X$

We want to find the minimal stabbing number Jordan curve containing $X$ in order to strengthen the lower bound on the cutting number from Lemma 3.1. In this section, we will talk about the minimal stabbing number Jordan curve through $X$.

**Definition 4.1 (Degree).** The degree $d(X)$ of $X$ is the smallest $d$ such that there exists a $d$-curve through $X$.

Using this definition, we can establish a better lower bound on the cutting number.

**Proposition 4.1.** *Given $X$, we have* $\mathrm{CUT}_K(X) \geq \max_{Y \subseteq X} \frac{|Y|}{K \cdot d(Y)}$.

*Proof.* By Lemma 3.1, for a $d(Y)$-curve through $Y$, we have $\mathrm{CUT}_K(X) \geq \mathrm{CUT}_K(Y) \geq \frac{|Y|}{K \cdot d(Y)}$. The result follows by considering this quantity over all $Y \subseteq X$. □

**Remark.** If $\max_{Y \subseteq X} \frac{|Y|}{K \cdot d(Y)} \gg \frac{N}{K^2}$, then $X$ cannot be separated equally.

This is the strongest condition we have discovered for determining $X$ cannot be separated equally. Since $d(X)$ features prominently in our result, we build our intuition for how $d(X)$ behaves. We start by exploring how large $d(X)$ can be, building off existing results for spanning trees. In 1988, Chazelle [2] proved the following theorem.

**Theorem 4.2 (Chazelle [2]).** *There exists a spanning tree through $X$ with* $\mathrm{stab}(X) = O(\sqrt{N})$.

The original proof for this theorem used a method called iterative reweighting, sequentially adding new edges by using a packing argument with respect to cuttings. Recently, new proofs have been found based on linear programming duality [9] and the Gutz-Katz polynomial partitioning technique [12]. With Theorem 4.2, we can bound $d(X)$.

**Lemma 4.3.** *For all $X$, we have $d(X) = O(\sqrt{N})$.*

*Proof.* By Theorem 4.2, there exists a spanning tree $T$ with $\mathrm{stab}(T) = O(\sqrt{N})$. We perform a pre-order traversal (depth-first search, traversing the left subtree completely before the right) of $T$ starting from an arbitrary root vertex in order to construct a closed curve $J$ around $T$ following this order. We keep $J$ at a fixed distance of arbitrarily small $\varepsilon > 0$ from $T$, except within $\varepsilon$ of vertices in $X$, where a pair of line segments connect $J$ to the vertex. An example is shown in Figure 4.1.

We claim $\mathrm{stab}(J) = O(\sqrt{N})$. Consider circles of radius $\varepsilon$ centered on each vertex. For $\varepsilon$ small enough, no line can intersect more than 2 of these circles. For every intersection of a line with $T$, there are at most 2 intersections with $J$ at a distance of $\varepsilon$ from $T$. This gives $O(\sqrt{N})$ intersections. If a line intersects $J$ at a distance of $< \varepsilon$ from $J$, it lies within one of the radius $\varepsilon$ circles. Since these circles only contain 2 line segments, there are at most 4 intersections at a distance of $< \varepsilon$ from $T$. This shows $\text{stab}(J) = O(\sqrt{N})$. $\square$

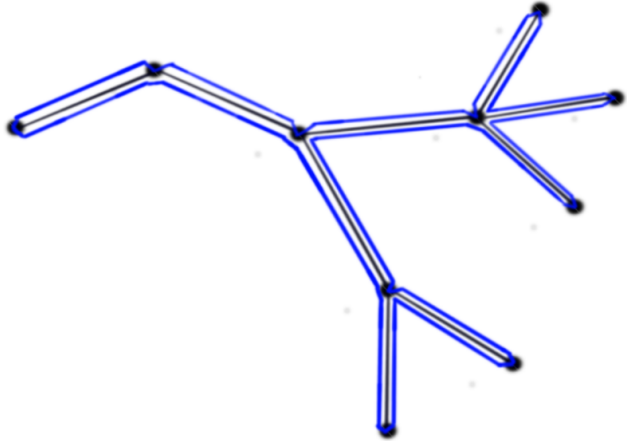

**Figure 4.1:** Example of $d$-curve (in blue) corresponding to spanning tree (in black) through $X$.

**Remark.** There exists $X$ with $d(X) = \Theta(\sqrt{N})$.

*Proof.* Consider an $X$ separated by $\mathcal{A}_K$ such that each cell has 1 point. Then there are $N$ cells so $K \sim \sqrt{N}$ and $\text{CUT}_K(X) = 1$. By Proposition 4.1, $\text{CUT}_K(X) = 1 \gtrsim \frac{N}{d(X) \cdot K} \sim \frac{\sqrt{N}}{d(X)}$. This implies $d(X) \gtrsim \sqrt{N}$ and by Lemma 4.3, the result follows. $\square$

We now have a worst-case optimal bound on $d(X)$. However, we can potentially bound the degree of a set better if we know the degree of certain subsets. If we have two linearly separable sets of points, then a line also separates their convex hulls. A polygon with vertices in $X$ is contained within the convex hull of $X$. Thus, in order to talk about the degree of the union of two linearly separable sets of points, we consider polygonal curves through $X$.

**Lemma 4.4.** *There exists a polygonal $d(X)$-curve through $X$ with vertices in $X$.*

*Proof.* Consider a $d(X)$-curve containing $X$ with points connected in the order $x_1, x_2, ..., x_n$. Construct a curve $J$ consisting of $\overline{x_1 x_2}, \overline{x_2 x_3}, ..., \overline{x_{n-1} x_n}, \overline{x_n x_1}$. Any line $\ell$ that intersects $\overline{x_i x_j}$ also intersects any arc joining $\overline{x_i x_j}$ at least once. Thus $\text{stab}(J) \leq \text{stab}(Z)$.

This curve $J$ is not guaranteed to be a polygon as it may intersect itself. If $\overline{x_i x_{i+1}}$ intersects $\overline{x_j x_{j+1}}$, we can delete these two segments of $J$ and add $\overline{x_i x_j}$ and $\overline{x_{i+1} x_{j+1}}$ to $J$. This removes the intersection and does not increase $\text{stab}(J)$. It also decreases the total length of the curve, so it follows that repeating this process eventually leads to a Jordan curve, giving us a polygonal $d(X)$-curve $J$ containing $X$ with vertices in $X$. □

Now we can bound the degree of the union of two linearly separable sets of points.

**Lemma 4.5.** *For linearly separable sets of points $A$ and $B$, we have $d(A \cup B) \leq d(A) + d(B) + 2$.*

*Proof.* By Lemma 4.4, we can draw a polygonal $d(A)$-curve $P$ through $A$ and a polygonal $d(B)$-curve $Q$ through $B$ such that there exists a line separating $P$ and $Q$. Add a line segment $R$ such that $|R \cap P| = |R \cap Q| = 1$ and $|R \cap A| = |R \cap B| = 0$. Then we can add a parallel line $S$ at an arbitrarily small distance from $R$ such that the subtended portions of $P$ and $Q$ contain no points in $A$ or $B$. Deleting these subtended portions gives us a new curve containing $A$ and $B$. An example is shown in Figure 4.2. Since we added 2 line segments, the final degree is at most $d(A) + d(B) + 2$. □

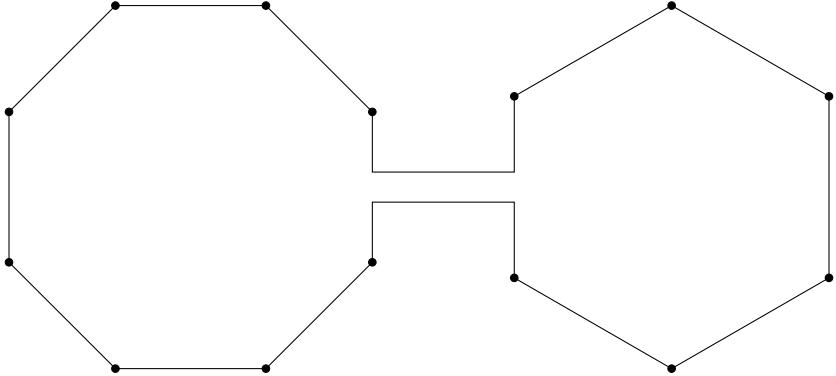

**Figure 4.2:** Example of $A$ and $B$ with $d(A) = d(B) = 2$ and $d(A \cup B) \leq 6$

**Theorem 4.6.** *Let $A_1, A_2, ..., A_m$ be pairwise linearly separable sets of points lying on corresponding $d(A_i)$-curves $B_1, B_2, ..., B_m$ for sufficiently large $m$. Then $d(\bigcup_i A_i) \leq \text{stab}(\bigcup_i B_i) + O(\sqrt{m})$.*

*Proof.* Place a point inside each $B_i$. By Theorem 4.2, there is a spanning tree $T$ through these points with stabbing number $O(\sqrt{m})$. If $T$ intersects points in $A_i$, we can perturb the original set of points inside $B_i$. For every $B_i$ and $B_j$ that $T$ connects, we can apply Lemma 4.5. The caveat is that the union of the line segments used to join the $B_i$ will have stabbing number $O(\sqrt{m})$. Any line in the plane can intersect $\bigcup_i B_i$ at most $\text{stab}(\bigcup_i B_i)$ times and the added line segments $O(\sqrt{m})$ times. Thus, we get $d(\bigcup_i A_i) \leq \text{stab}(\bigcup_i B_i) + O(\sqrt{m})$. □

# 5 Generalizing the Erdős-Szekeres Conjecture

The Happy Ending Problem has a rich history as a motivating problem in the origins of Ramsey theory. Even its name has a great story: Erdős coined it because collaboration over the problem between George Szekeres and Esther Klein led to their marriage. The problem states that every set of 5 points in general position has a subset of 4 points that form a convex quadrilateral. In our notation, this is equivalent to saying that if $N \geq 5$, then $\text{CON}(X) \geq 4$. In 1935, Erdős and Szekeres [5] generalized this statement.

**Definition 5.1 (F(n)).** We let $F(n)$ be the minimum number of points such that there always exists a subset of $n$ points that form a convex $n$-gon.

**Conjecture 5.1 (Erdős-Szekeres [5]).** *For $n \geq 3$, we have $F(n) = 2^{n-2} + 1$.*

Erdős and Szekeres further showed [6] that $F(n) \geq 2^{n-2} + 1$. This gives us a lower bound on the convex number, namely that $\text{CON}(X) \geq \log_2(4N - 4)$. However, the best upper bound they could muster was $F(n) \leq \binom{2n-4}{n-2} + 1 \sim \frac{4^n}{\sqrt{n}}$ by Ramsey's Theorem [5]. For 81 years, small improvements [14] were made to the bound, but the order of magnitude remained unchanged. In April 2016, Suk [15] showed the following and made significant progress toward proving the conjecture.

**Theorem 5.2 (Suk [15]).** *We have $F(n) \leq 2^{n+o(n)}$.*

However, Theorem 5.2 only tells us about the existence of convex curves. We can generalize this result to $d$-curves.

**Definition 5.2 (F(n, d)).** We let $F(n,d)$ be the minimum number of points such that there always exists a subset of $n$ points that lie on a $d$-curve.

**Theorem 5.3.** *For fixed $d$ and sufficiently large $n$, $F(n,d) \leq 2^{c_d n/d + o(n)}$ with $c_d = 1 + O(\frac{1}{\sqrt{d}})$.*

*Proof.* Let $h = \frac{d}{2} - C\sqrt{d}$ for some constant $C$ to be chosen later. Then as $d \to \infty$, we have $h \to \infty$.

Take $X$ with at least $(2h+2)(2^{\lceil n/2h \rceil + o(\lceil n/2h \rceil)})$ points. Apply a partial cutting with $h+1$ lines. At least $2h$ of the cells formed will contain at least $2^{\lceil n/2h \rceil + o(\lceil n/2h \rceil)}$ points. By Theorem 5.2, each of these $2h$ cells contains a convex curve $B_i$ with $\lceil n/2h \rceil$ points. Let $A$ be the union of all the points on these convex curves. Then as $d \to \infty$, by Theorem 4.6, we have $d(A) \leq \text{stab}(\bigcup_i B_i) + O(\sqrt{h})$.

We bound $\text{stab}(\bigcup_i B_i)$. Any line can intersect at most $h+1$ of the lines of the partial cutting. Hence, any line can only intersect $h+2$ cells and intersect $\bigcup_i B_i$ twice in each of them. This is a total of at most $2h+4$ intersections with the convex curves, giving $\text{stab}(\bigcup_i B_i) \leq 2h+4$.

Thus, we have $d(A) \leq 2h + O(\sqrt{h}) \leq d$. We pick a $C$ such that $h = \frac{d}{2} - C\sqrt{d}$ satisfies this condition as $d \to \infty$. This gives us:

$$\begin{align}
F(n,d) &\leq (2h+2)(2^{\lceil n/2h \rceil + o(\lceil n/2h \rceil)}) \\
&= (d - 2C\sqrt{d} + 2)(2^{\lceil n/(d-2C\sqrt{d}) \rceil + o(\lceil n/(d-2C\sqrt{d}) \rceil)}) \tag{1} \\
&= 2^{\frac{n}{d-2C\sqrt{d}} + o(n)} \tag{2} \\
&= 2^{\frac{n}{d}(\frac{d+2C\sqrt{d}}{d-4C^2}) + o(n)} \tag{3} \\
&= 2^{\frac{n}{d}(1 + O(\frac{1}{\sqrt{d}})) + o(n)} \tag{4}
\end{align}$$

We go from (1) to (2) because $n$ is sufficiently large ($n \gg d$), allowing the $o(n)$ term to absorb most of the terms involving $d$ and the ceiling functions. Simple algebraic manipulations take us from (2) to (3). Finally, as $d \to \infty$, we obtain (4), completing the proof. □

For small $d$, the following is also relevant.

**Corollary 5.3.1.** *We also have $c_d \leq d$.*

*Proof.* Since convex curves are also $d$-curves, we have

$$2^{c_d n/d + o(n)} = F(n,d) \leq F(n,2) = 2^{2n/2 + o(n)}$$

which gives $c_d \leq d$. □

Suk's result in Theorem 5.2 was that $c_2 = 2$. We showed that $c_d$ approaches 1 as $d$ goes to infinity. In the process, we drew together various ideas discussed in prior sections, such as partial cuttings and degrees of subsets. These coalesced in a natural generalization of a major problem in Ramsey theory, giving us insight into the existence of $d$-curves.

# 6 Conclusions and Future Work

## 6.1 Future research

Many related problems still remain. Obviously, we would like to identify all obstacles to equal separation. Experimentation using computers could be fruitful, but significant difficulties arise.

For example, given $X$, we would need to find a $d(X)$-curve through $X$ to explore the relationship between $d$-curves and equal separation.

**Problem 6.1.** *What is the complexity of computing a $d(X)$-curve through $X$?*

The NP-hardness of computing spanning trees of minimal stabbing number [7] would imply that Problem 6.1 is challenging too. If we were able to compute $d(X)$ quickly, we could generate random point sets and see how frequently $d(X)$ is large. As we explained in Section 1, being able to equally separate a set of points should not be a common occurrence. If large $d(X)$ allows for equal separation, we would expect $d(X)$ to be small more often. Obviously, to check if $X$ can be equally separated, the following problem is relevant.

**Problem 6.2.** *If $X$ can be separated equally by some $\mathcal{A}_K$, how do we find such an $\mathcal{A}_K$?*

It can be shown $O(N^2)$ lines check all possible separations of $X$ with one line. It follows that computing $\mathcal{A}_K$ should run in $O(N^{2K})$ time. This is too large for even relatively small $N$ and $K$.

Since it seems difficult to compute $d(X)$ or an $\mathcal{A}_K$ separating $X$ equally, we might want to show $d(X)$ is large without considering all of $X$. For example, to check if $X$ lies on a convex curve, one can check all subsets of 4 points. If some 4-point subset of $X$ does not lie on a convex curve, then $X$ does not lie on a convex curve and $d(X) > 2$. The proof is simple. If $X$ does not lie on a convex curve, then there is some point inside its convex hull. This point lies in a triangle formed by 3 other points of $X$, so these 4 points do not form a convex quadrilateral. However, this relies on the fact that a convex hull can be triangulated, which does not generalize easily to $d$-curves for $d > 2$.

**Problem 6.3.** *Given some $d$, does there exist some $k$ such that if all subsets of $k$ points in $X$ lie on a $d$-curve then $X$ lies on a $d$-curve?*

An answer to Problem 6.3 might let us show that $d(X)$ is large without explicitly computing curves through $X$, allowing us to collect relevant information from computer simulations without solving Problem 6.1.

The original motivating question of this paper was to determine the obstacles to equal separation. Proposition 4.1 gave the best known criterion for when $X$ cannot be equally separated. Its converse would establish low degree as the only obstacle to equal separation.

**Conjecture 6.4.** *If $X$ cannot be separated equally, then $\max\limits_{Y \subseteq X} \frac{|Y|}{K \cdot d(Y)} \gg \frac{N}{K^2}$.*

The author is unaware of a set with subsets of large degree that does not separate equally. As mentioned earlier, finding a counterexample may be aided by computational experiments. Together, these open problems and conjectures leave many paths for building on results from this paper.

## 6.2 Implications and relevance

Our work improved on the known obstacles to equal separation. We established the criteria that if $\max_{Y \subseteq X} \frac{|Y|}{K \cdot d(Y)} \gg \frac{N}{K^2}$, then $X$ cannot be separated equally. In the process, we explored a variety of related ideas. We developed a method to separate $X$ with $K$ lines into cells with at most $\lceil \frac{N}{2K} \rceil$ points. We showed that there exists a polygon through $X$ with vertices in $X$ and stabbing number $O(\sqrt{N})$. Finally, we generalized the Erdős-Szekeres Conjecture to show that for sufficiently large $n$, we have $F(n,d) \leq 2^{c_d n/d + o(n)}$ with $c_d = 1 + O(\frac{1}{\sqrt{d}})$.

Guth [8] mentions generalizing the Szemerédi-Trotter theorem to $\mathbb{R}^n$ for $n \geq 3$. Indeed, we can equally distribute $(n-1)$-dimensional planes in $\mathbb{R}^n$ among $\sim K^n$ cells with an arrangement of $\sim K$ hyperplanes. The argument is very similar to the cell decomposition method for distributing lines in $\mathbb{R}^2$. However, he notes distributing $m$-dimensional objects in $\mathbb{R}^n$ for $m < n-1$ breaks down in higher dimensions just as it does for $m = 0$ and $n = 2$. The conclusion is problems dealing with objects of codimension greater than 1 are far more challenging. Our approach to the problem involved dealing with Jordan curves. These are objects of codimension 1, which makes equal distribution a more tractable problem and fits with these observations by Guth.

As mentioned in Section 1, stabbing numbers are a well-studied topic in computational geometry. The bulk of relevant research has been focused on spanning trees, triangulations, matchings, and relevant algorithmic development. We showed the significance of minimal stabbing number Jordan curves to a problem in combinatorial geometry. Our work might be used in other problems utilizing low stabbing number objects, such as topics cited in Section 1, as well as the generalizations discussed in the previous paragraph. Meanwhile, arrangements of hyperplanes are essential in various subjects, including combinatorial geometry, Lie algebras, and Coxeter groups [11]. Our research discussed fundamental questions about how arrangements of lines behave, presenting new questions for future research.

# 7    Acknowledgements

The author would like to thank his mentor Borys Kadets (MIT) for invaluable insight and guidance. He would like to thank Dr. Larry Guth (MIT) for suggesting this project. He would also like to thank Dr. Tanya Khovanova (MIT) and Dr. Slava Gerovitch (MIT) for helpful advice and suggestions. Finally, he would like to thank the MIT PRIMES-USA program for providing this research opportunity.